\documentclass[12pt, reqno]{amsart}
\usepackage{amsmath}
\usepackage{amsxtra}
\usepackage{amscd}
\usepackage{amsthm}
\usepackage{amsfonts}
\usepackage{amssymb}
\usepackage{eucal}

\newtheorem{theorem}{Theorem}[section]

\newtheorem{cor}[theorem]{Corollary}
\newtheorem{lem}[theorem]{Lemma}
\newtheorem{prop}[theorem]{Proposition}

\theoremstyle{definition}

\theoremstyle{remark}
\newtheorem{rem}[theorem]{Remark}
\theoremstyle{remark}

\numberwithin{equation}{section}

\newcommand{\nc}{\newcommand}
\nc{\on}{\operatorname}
\nc{\ch}{\mbox{ch}}
\nc{\Z}{{\mathbb Z}}
\nc{\C}{{\mathbb C}}
\nc{\pone}{{\mathbb C}{\mathbb P}^1}
\nc{\pa}{\partial}
\nc{\F}{{\mathcal F}}
\nc{\arr}{\rightarrow}
\nc{\larr}{\longrightarrow}
\nc{\al}{\alpha}
\nc{\ri}{\rangle}
\nc{\lef}{\langle}
\nc{\W}{{\mathcal W}}
\nc{\la}{\lambda}
\nc{\ep}{\epsilon}
\nc{\su}{\widehat{{\mathfrak sl}}_2}
\nc{\sw}{{\mathfrak s}{\mathfrak l}}
\nc{\g}{{\mathfrak g}}
\nc{\h}{{\mathfrak h}}
\nc{\n}{{\mathfrak n}}
\nc{\N}{\widehat{\n}}
\nc{\G}{\widehat{\g}}
\nc{\De}{\Delta_+}
\nc{\gt}{\widetilde{\g}}
\nc{\Ga}{\Gamma}
\nc{\one}{{\mathbf 1}}
\nc{\z}{{\mathfrak Z}}
\nc{\zz}{{\mathcal Z}}
\nc{\Hh}{{\mathcal H}_\beta}
\nc{\qp}{q^{\frac{k}{2}}}
\nc{\qm}{q^{-\frac{k}{2}}}
\nc{\La}{\Lambda}
\nc{\wt}{\widetilde}
\nc{\qn}{\frac{[m]_q^2}{[2m]_q}}
\nc{\cri}{_{\on{cr}}}
\nc{\kk}{h^\vee}
\nc{\sun}{\widehat{\sw}_N}
\nc{\hh}{\widehat{\mathfrak h}}
\nc{\HH}{{\mathcal H}_{q,t}}
\nc{\ca}{\wt{{\mathcal A}}_{h,k}(\sw_2)}
\nc{\gl}{\widehat{{\mathfrak g}{\mathfrak l}}_2}
\nc{\el}{\ell}
\nc{\s}{{\mathbf s}}
\nc{\bi}{\bibitem}
\nc{\om}{\omega}
\nc{\WW}{\W_\beta}
\nc{\scr}{{\mathbf S}}
\nc{\ab}{{\mathbf a}}
\nc{\rr}{r}
\nc{\ol}{\overline}
\nc{\con}{qt^{-1} + q^{-1}t}
\nc{\den}{q^{\el-1} t^{-\el+1}+ q^{-\el+1} t^{\el-1}}
\nc{\ds}{\displaystyle}
\nc{\B}{B}
\nc{\A}{{\mathbb A}}
\nc{\GG}{{\mathcal G}}
\nc{\UU}{{\mathcal U}}
\nc{\MM}{{\mathcal M}}
\nc{\CC}{{\mathcal C}}
\nc{\GL}{{}^L G}
\nc{\dzz}{\frac{dz}{z}}
\nc{\Res}{\on{Res}}
\nc{\rep}{{\mathcal R}ep \;}
\nc{\uqg}{U_q \G}
\nc{\uqgg}{U_q \g}
\nc{\Fq}{{\mathbb F}_q}
\nc{\sgn}{\on{sgn}}
\nc{\stimes}{\ltimes}
\nc{\K}{\hat{\mathcal K}}
\nc{\Ql}{\ol{\mathbb Q}_\ell}

\nc{\ga}{\gamma}
\nc{\PL}{{}^L P}
\nc{\E}{\mc E}
\nc{\mc}{\mathcal}
\nc{\mbf}{\mathbf}
\nc{\bb}{{\mathfrak b}}
\nc{\OO}{{\mc O}}
\nc{\Po}{{\mc P}}
\nc{\V}{{\mc V}}
\nc{\yy}{{\mc Y}}
\nc{\M}{\mathcal M}
\nc{\Coh}{{{\mathcal C}oh}}
\nc{\Cohn}{\Coh_n}
\nc{\f}{{\mathcal F}}
\nc{\si}{_E}
\nc{\Gaf}{{\mathbb G}_{a,\Fq}}
\nc{\KK}{{\mathfrak k}}
\nc{\pr}{\partial}
\def\gl{\frak{gl}}

\newcommand{\NNN}{^{( N )}}
\newcommand{\MMM}{^{( M )}}
\def\Vb{L_\bullet}

\newcommand{\bean}{\begin{eqnarray}}
\newcommand{\eean}{\end{eqnarray}}
\newcommand{\be}{\begin{displaymath}}
\newcommand{\ee}{\end{displaymath}}
\newcommand{\bea}{\begin{eqnarray*}}   
\newcommand{\eea}{\end{eqnarray*}}
\newcommand{\bs}{\boldsymbol}
\newcommand{\Ref}[1]{{$($\ref{#1}$)$}}

\begin{document}
\title[A generalization of the Capelli identity]
{A generalization of the Capelli identity}

\author[E. Mukhin, V. Tarasov, and
A. Varchenko] {E. Mukhin, V. Tarasov, and A. Varchenko}
\thanks{Research of E.M. is supported in part by NSF grant
DMS-0601005.  Research of A.V. is supported in part by NSF grant
DMS-0555327} \address{E.M.: Department of Mathematical Sciences,
Indiana University -- Purdue University Indianapolis, 402 North
Blackford St, Indianapolis, IN 46202-3216, USA,\newline
mukhin@math.iupui.edu} \address{V.T.: Department of Mathematical
Sciences, Indiana University -- Purdue University Indianapolis, 402
North Blackford St, Indianapolis, IN 46202-3216, USA,\newline
vtarasov@math.iupui.edu, \ and\newline \hglue\parindent St.~Petersburg
Branch of Steklov Mathematical Institute, Fontanka 27,\newline
St.~Petersburg, 191023, Russia, vt@pdmi.ras.ru} \address{A.V.:
Department of Mathematics, University of North Carolina at Chapel
Hill,\newline Chapel Hill, NC 27599-3250, USA, anv@email.unc.edu}

\setcounter{footnote}{0}\renewcommand{\thefootnote}{\arabic{footnote}}

\begin{abstract}
We prove a generalization of the Capelli identity.
As an application we obtain an isomorphism of the Bethe subalgebras
actions under the $(\gl_N,\gl_M)$ duality.
\end{abstract}
\maketitle

\section{Introduction}
Let $\mc A$ be an associative algebra over complex numbers. 
Let $A=(a_{ij})_{i,j=1}^n$ be an $n\times n$ 
matrix with entries in $\mc A$.
The {\it row determinant of $A$} is defined by the formula:
\be
\on{rdet}(A):=
\sum_{\sigma \in S_n} \sgn(\sigma)a_{1\sigma_1}\dots a_{n\sigma_n}.
\ee

Let $x_{ij}$, $i,j=1,\dots,M$, be commuting variables. 
Let $\pr_{ij}=\pr/\pr x_{ij}$, 
\bean
\label{left act}
E_{ij}=\sum_{a=1}^M x_{ia}\pr_{ja}.
\eean
Let $X=(x_{ij})_{i,j=1}^M$ and  
$D=(\pr_{ij})_{i,j=1}^M$ be $M\times M$ matrices.

The classical Capelli  identity \cite{C1} asserts the following equality 
of differential operators:
\bean\label{capelli}
\on{rdet} \Big(E_{ji}+(M-i)\delta_{ij}\Big)_{i,j=1}^M=\det(X) \det (D).
\eean
This identity is a ``quantization'' of the identity 
$$
\det(AB) =\det(A) \det(B)
$$ 
for any matrices $A,B$ with commuting entries.

\medskip

The Capelli identity has the following meaning in the representation
theory.  Let $\C[X]$ be the algebra of complex polynomials in
variables $x_{ij}$.  There are two natural actions of the
Lie algebra
$\gl_M$ on $\C[X]$.  The first action is given by operators from
 \Ref{left act} and
the second action is given by operators
$\widetilde  E_{ij}=\sum_{a=1}^M x_{ai}\pr_{aj}.$  
The two actions commute
and the corresponding $\gl_M\oplus \gl_M$ action is multiplicity free.

It is not difficult to see that the right hand side of \Ref{capelli}, 
considered as a differential operator on $\C[X]$, commutes with
both actions of $\gl_M$ and therefore lies in the image of the center
of the universal enveloping algebra $U\gl_M$ with respect to the first action.
Then the left hand side of the Capelli
identity  expresses the corresponding central element in terms of 
$U\gl_M$ generators.

\medskip

Many generalizations of the Capelli identity 
are known. 
One group of 
generalizations  considers 
other elements of the center of $U\gl_M$, 
called quantum immanants,
and then expresses them in terms of $\gl_M$ generators, see \cite{C2}, 
\cite{N1},\cite{O}. 
Another group of generalizations considers other 
pairs of Lie algebras in place of $(\gl_M,\gl_M)$, 
 e.g. $(\gl_M,\gl_N)$, $(\frak{sp}_{2M},\gl_2)$, 
$(\frak{sp}_{2M},\frak{so}_N)$, etc,  see \cite{MN}, \cite{HU}. 
The third group 
of generalizations produces 
identities corresponding not to pairs of Lie algebras,
but to pairs of
quantum groups \cite{NUW} or superalgebras \cite{N2}.

In this paper we prove a generalization of the Capelli identity which 
seemingly does not fit the above classification.

\medskip

Let $\bs z=(z_1,\dots,z_N)$, $\bs \la=(\la_1,\dots,\la_M)$ be sequences
of complex numbers.
Let $Z=(z_i\delta_{ij})_{ij=1}^N$, $\La=(\la_i\delta_{ij})_{ij=1}^M$ be 
the corresponding diagonal matrices. Let $X$ and $D$ be the 
$M\times N$ matrices with  entries
$x_{ij}$ and $\pr_{ij}$, $i=1,\dots, M$, $j=1,\dots,N$, respectively.
Let $\C[X]$ be the algebra of complex polynomials in
variables $x_{ij}$, $i=1,\dots, M$, $j=1,\dots,N$.
Let $E_{ij}^{(a)}=x_{ia}\pr_{ja}$,
where $i,j=1,\dots,M,\,a=1,\dots,N$.

In this paper we prove that
\bean
\label{new}
\prod_{a=1}^N(u-z_a)\ \on{rdet}\Big((\pr_u-\la_i)\delta_{ij}-
\sum_{a=1}^N\frac{E_{ji}^{(a)}}{u-z_a}\Big)_{i,j=1}^M=
\on{rdet} \left( 
\begin{matrix}
u-Z & X^t \\
D & \pr_u-\La
\end{matrix}
\right)\ .
\eean
The left hand side of \Ref{new} is an
 $M\times M$ matrix while the right hand side 
is an $(M+N)\times (M+N)$ matrix.

Identity \Ref{new} is a ``quantization'' of the identity 
\be
\det \left(
\begin{matrix}
A & B \\
C & D
\end{matrix}
\right)=\det(A)\  \det (D-CA^{-1}B)
\ee 
which holds for any matrices $A,B,C,D$ with commuting entries,  
for the case when $A$ and $D$ are  diagonal matrices.

By setting all  
$z_i,\la_j$ and $u$ to zero, and $N=M$ in \Ref{new}, 
we obtain the classical Capelli identity
\Ref{capelli}, see Section \ref{cap}.

\medskip

Our proof of \Ref{new} is combinatorial and reduces 
to the case of $2\times 2$ matrices. 
In particular, it gives  a 
proof of the classical
Capelli identity, which may be new.

\medskip

We invented  identity \Ref{new}  to prove Theorem \ref{dual} below,  and
   Theorem  \ref{dual} in its turn was motivated by results of \cite{MTV2}.
   In Theorem \ref{dual} we compare  actions of two Bethe subalgebras.

Namely, consider
$\C[X]$ as a tensor product of evaluation modules over the current Lie
algebras $\gl_M[t]$
 and $\gl_N[t]$ with evaluation
parameters $\bs z$ and $\bs \la$, respectively. 
The action of the algebra 
$\gl_M[t]$ on $\C[X]$ is given by the formula
\be
E_{ij}\otimes t^n=\sum_{a=1}^N x_{ia}\pr_{ja}z_i^n,
\ee
and the  action of the  algebra 
$\gl_N[t]$ on $\C[X]$ is given by the formula
\be
E_{ij}\otimes t^n=\sum_{a=1}^M x_{ai}\pr_{aj}\la_i^n.
\ee
In contrast to the previous situation, these two actions do not commute.

The algebra $U\gl_M[t]$ has a family of commutative subalgebras $\mc
G(M,\bs \la)$ depending on parameters $\bs \la$ and called the Bethe
subalgebras.  For a given $\bs \la$, the Bethe subalgebra $\mc G(M,\bs
\la)$ is generated by the coefficients of the expansion of the
expression \bean\label{Gaudin} \on{rdet}\Big((\pr_u-\la_i)\delta_{ij}-
\sum_{a=1}^N \sum_{s=1}^\infty (E_{ji}^{(a)}\otimes
t^s)u^{-s-1}\Big)_{i,j=1}^M \eean with respect to powers of $u$ and
$\pr_u$, cf. Section 3. For different versions of definitions of Bethe
subalgebras and relations between them, see \cite{FFR}, \cite{T},
 \cite{R}, \cite{MTV}.

Similarly, there is a family of Bethe subalgebras $\mc G(N,\bs z)$ in 
$U\gl_N[t]$ depending on parameters $\bs z$.

For fixed $\bs\la$ and $\bs z$, consider the action of the Bethe subalgebras
$\mc G(M,\bs \la)$ and $\mc G(N,\bs z)$ on $\C[X]$ as defined above.
In Theorem \ref{dual} we show that  the actions of the Bethe subalgebras
on  $\C[X]$ induce the same subalgebras of endomorphisms of $\C[X]$.

\medskip

The  paper is organized as follows. 
In Section \ref{identity sec} we describe  and prove formal Capelli-type
identities and in
Section \ref{bethe sec} we discuss the relations of the identities
 to the Bethe subalgebras.

\section{Identities}\label{identity sec}
\subsection{The main identity}
We work over the field of complex numbers, however 
all results of this paper hold over any field of 
characteristic zero.

Let $\mc A$ be an associative algebra.
Let $A=(a_{ij})_{i,j=1}^n$ be an $n\times n$ 
matrix with entries in $\mc A$.
Define the {\it row determinant of $A$} by the formula:
\be
\on{rdet}(A):=
\sum_{\sigma \in S_n} \sgn(\sigma)a_{1\sigma_1}\dots a_{n\sigma_n},
\ee 
where $S_n$ is the symmetric group on $n$ elements.

Fix two natural numbers $M$ and $N$ and a complex number $h\in\C$. Consider
noncommuting variables $u, p_u, x_{ij}, p_{ij}$, where 
$i=1,\dots, M$, $j=1,\dots ,N$, such that 
the commutator of two variables equals zero except 
\be
[p_u,u]=h, \qquad [p_{ij},x_{ij}]=h,
\ee
$i=1,\dots, M$, $j=1,\dots, N$.

Let $X,P$ be two $M\times N$ matrices given by 
\be
X\,:\ =\ (x_{ij})_{i=1,\dots,M}^{j=1,\dots,N}\ , 
\qquad 
P\,:\ =\ (p_{ij})_{i=1,\dots,M}^{j=1,\dots,N}\ .
\ee

Let $\mc A_h^{(MN)}$ be the associative algebra 
whose elements are polynomials in
$p_u, x_{ij}, p_{ij}$, $i=1,\dots, M$, $j=1,\dots, N,$ with
coefficients that are rational functions in $u$. 

Let $\mc A^{(MN)}$ be the associative algebra of 
linear differential operators in $u, x_{ij}$, $i=1,\dots,
M$, $j=1,\dots, N$, with coefficients in $\C(u)\otimes \C[X]$.

We often drop the dependence on $M,N$ and write  $\mc A_h$, $\mc A$
for $\mc A^{(MN)}_h$ and $\mc A^{(MN)}$, respectively.

For $h\neq 0$, we have the isomorphism of algebras
\bean\label{isom}
\iota_h\ :\ \mc A_h &\to& \mc A\ , 
\\ 
u,\,x_{ij}& \mapsto & u,\,x_{ij}\ , 
\notag
\\
p_u, \,p_{ij} & \mapsto & h \frac{\pr}{\pr u},\, h\frac{\pr}{\pr x_{ij}}
\notag \ .
\eean

Fix two sequences of complex numbers $\bs z=(z_1,\dots,z_N)$ and 
$\bs \la=(\la_1,\dots,\la_M)$.

Define the $M\times M$ matrix 
$G_h\,=\,G_h(M,N,u,p_u,\bs z,\bs \la,X,P)$ by the formula
\bean\label{G}
G_h\,:\ =\  
\Big(\,(p_u-\la_i)\,\delta_{ij}\ -\ 
\sum_{a=1}^N\,\frac{x_{ja}p_{ia}}{u-z_a}\,\Big)_{i,j=1}^M\ .
\eean

\begin{theorem}\label{main}
We have
\be
\prod_{a=1}^N(u-z_a)\ \on{rdet}(G_h)=\hspace{-10pt}
\sum_{A,B,  |A|=|B|} (-1)^{|A|}\prod_{a\not\in B}(u-z_a)
\prod_{b\not \in A}(p_u-\la_b) \ 
\det(x_{ab})_{a\in A}^{b\in B} \ 
\det(p_{ab})_{a\in A}^{b\in B},
\ee
where the sum is over all pairs of subsets $A\subset\{1,\dots,M\}$, 
$B\subset\{1,\dots,N\}$ such that $A$ and $B$ have 
the same cardinality, $|A|=|B|$. Here the sets $A,B$ 
inherit the natural ordering from the sets $\{1,\dots,M\}$, $\{1,\dots,N\}$. 
This ordering determines the determinants in the formula.
\end{theorem}

Theorem \ref{main} is proved in Section \ref{proof}.

\subsection{A presentation as a row determinant of size $M+N$}
Theorem \ref{main} implies 
that the row determinant of $G$ can be written as the
row determinant of a matrix of size $M+N$.

Namely, let  $Z$ be
the diagonal $N\times N$ matrix with diagonal entries $z_1,\dots,z_N$.
Let $\La$ be the diagonal $M\times M$ matrix with diagonal entries
$\la_1,\dots,\la_M$:
\be
Z\,:\ =\ (\,z_i\delta_{ij}\,)_{i,j=1}^N\ , 
\qquad 
\La\,:\ =\ (\,\la_i\delta_{ij}\,)_{i,j=1}^M\ .
\ee 
\begin{cor}\label{m+n}
We have 
\be
\prod_{a=1}^N(u-z_a)\ \on{rdet} \ G=\on{rdet} \left( 
\begin{matrix}
u-Z & X^t \\
P & p_u-\La
\end{matrix}
\right),
\ee
where $X^t$ denotes the transpose of the matrix $X$.
\end{cor}
\begin{proof}
Denote
\be
W\,:\ =\  \left( \begin{matrix}
u-Z & X^t \\
P &p_u-\La
\end{matrix}
\right), 
\ee 
The entries of the first $N$ rows of $W$ commute. The
entries of the last $M$ rows of $W$ also commute.  
Write the Laplace
decomposition of $\on{rdet}(W)$ with respect to the first $N$ rows.
Each term in this decomposition corresponds to a choice of $N$ columns
in the $N\times (N+M)$ matrix $(u-Z,X^T)$. We label such a choice by a
pair of subsets $A\subset\{1,\dots,M\}$ and $B\subset\{1,\dots,N\}$ of the same
cardinality.
Namely, the elements of $A$ correspond to the chosen columns in $X^T$ and the 
elements of the complement
to $B$ correspond to the chosen columns in $u-Z$. Then the term in the Laplace
decomposition corresponding to $A$ and $B$ is exactly the term labeled
by $A$ and $B$ in the right hand side of the formula in Theorem
\ref{main}. Therefore, the corollary follows from Theorem \ref{main}.
\end{proof}
Let $A, B, C, D$ be any matrices with commuting entries 
of sizes $N\times N, N\times M, M\times N$ and $M\times M$, respectively. 
Let $A$ be invertible. Then we have the equality of matrices of sizes 
$(M+N)\times(M+N)$:
\be
\left(\begin{matrix}
A & B \\
C & D
\end{matrix}\right)=
\left(\begin{matrix}
A & 0 \\
C & D-CA^{-1}B
\end{matrix}\right)
\left(\begin{matrix}
1 & A^{-1}B \\
0 & 1
\end{matrix}\right)
\ee
and therefore
\bean\label{gauss}
\det\left( 
\begin{matrix}
A & B \\
C & D
\end{matrix}
\right)\ =\ \det(A)\,  \det (D-CA^{-1}B)\ .
\eean
The identity of
Corollary \ref{m+n} for $h=0$ turns into identity \Ref{gauss} 
with diagonal matrices $A$ and $D$.
 Therefore, the identity of Corollary \ref{m+n}
may be thought of as a ``quantization'' 
of identity \Ref{gauss} with diagonal  $A$ and $D$.


\subsection{A relation between determinants of sizes $M$ and $N$}
Introduce new variables $v,p_v$ such that $[p_v,v]=h$.

Let $\bar{\mc A}_h$ be the associative algebra 
whose elements are polynomials in
$p_u, p_v, x_{ij}, p_{ij}$, $i=1,\dots, M$, $j=1,\dots, N,$ with
coefficients  in $\C(u)\otimes \C(v)$.

Let $e:\bar{\mc A}_h \to \bar{\mc A}_h$ be the unique linear map 
which
is the identity map on the subalgebra of $\bar{\mc A}_h$ generated by all
monomials which do not contain $p_u$ and $p_v$ and which satisfy 
\be
e(ap_u)=e(a)v, \qquad e(ap_v)= e(a)u, 
\ee 
for any $a\in\bar{\mc A}_h$.

Let $\bar{\mc A}$ be the associative algebra of
linear differential operators in $u, v, x_{ij}$, $i=1,\dots,
M$, $j=1,\dots, N$, with coefficients in $\C(u)\otimes \C(v)\otimes \C[x_{ij}]$.
Then for $h\neq 0$, we have the isomorphism of algebras extending 
the isomorphism \Ref{isom}:
\bea
\bar{\iota}_h:\ \bar{\mc A_h} &\to& \bar{\mc A}, \\ 
u,v, x_{ij}&\mapsto& u,v,x_{ij}, \notag\\
p_u,p_v, p_{ij}&\mapsto& h\frac{\pr}{\pr u}, h\frac{\pr}{\pr v}, 
h\frac{\pr}{\pr x_{ij}}\notag \ .
\eea

For $a\in \bar{\mc A}$ and a function $f(u,v)$ let $a\cdot f(u,v)$
denotes the function obtained by the action of $a$ considered as a 
differential operator in $u$ and $v$ on the function $f(u,v)$.

We have 
\be 
\bar\iota_h(e(a))=\exp(-uv/h)\bar\iota_h(a)\cdot\exp(uv/h)
\ee 
for any $a\in\bar A_h$ such that $a$ does not depend on either
$p_u$ or $p_v$.

Define the $N\times N$ matrix $H_h=H_h(M,N,v,p_v,\bs z,\bs \la,X,P)$ by
\bean\label{H}
H_h:=\Big((p_v-z_i)\delta_{ij}-
\sum_{b=1}^M\frac{x_{bj}p_{bi}}{v-\la_b}\Big)_{i,j=1}^N,
\eean
cf. formula \Ref{G}.
\begin{cor}\label{duality rel}
We have
\be
e\Big(\prod_{a=1}^N(u-z_a) \on{rdet}(G_h)\Big)= 
e\Big(\prod_{b=1}^M(v-\la_b)\on{rdet}(H_h)\Big).
\ee
\end{cor}
\begin{proof}
Write the dependence on parameters of the matrix $G$:
$G_h=G_h(M,N,u,p_u,\bs z,\bs \la,X,P)$. Then
\be
H_h=G_h(N,M,v,p_v,\bs \la,\bs z,X^T,P^T).
\ee
The corollary now follows from Theorem \ref{main}. 
\end{proof}

\subsection{A relation to the Capelli identity}\label{cap}
In this section we show how to deduce the Capelli identity from Theorem 
\ref{main}.

Let $s$ be a complex number.
Let $\al_s:{\mc A}_h \to {\mc A}_h$ be the unique linear map which
is the identity map on the subalgebra of ${\mc A}_h$ generated by all
monomials which do not contain $p_u$, and which satisfies 
\be
\al_s(aup_u)= s\al_s(a)
\ee 
for any $a\in\bar{\mc A}_h$.

We have
\be
\bar\iota_h(\al_s(a))=u^{-s/h}\bar\iota_h(a)\cdot u^{s/h}
\ee
for any $a\in\bar A_h$.

Consider the case $z_1=\dots=z_N=0$ and $\la_1=\dots=\la_M=0$ in Theorem 
\ref{main}.

Then it is easy to see that
the row determinant $\on{rdet}(G)$ can be rewritten in the following form
\be
u^M \on{rdet}(G_h)=
\on{rdet}
\Big(h(up_u-M+i)\delta_{ij}-\sum_{a=1}^N x_{ja}p_{ia}\Big)_{i,j=1}^M.
\ee 
Applying the map $\al_s$, we get
\be
\al_s(u^M \on{rdet}(G_h))=
\on{rdet}
\Big(h(s-M+i)\delta_{ij}-\sum_{a=1}^N x_{ja}p_{ia}\Big)_{i,j=1}^M.
\ee
Therefore applying Theorem \ref{main} we obtain the identity
\be
\on{rdet}\Big(h(s-M+i)\delta_{ij}-\sum_{a=1}^N x_{ja}p_{ia}\Big)_{i,j=1}^M=
\hspace{-15pt}
\sum_{A,B,  |A|=|B|} \hspace{-5pt}(-1)^{|A|}\hspace{-8pt}
\prod_{b=0}^{M-|A|-1}\hspace{-8pt}(s-bh) \ 
\det(x_{ab})_{a\in A}^{b\in B} \ 
\det(p_{ab})_{a\in A}^{b\in B}.
\ee
In particular, if $M=N$, and $s=0$, 
we obtain the famous Capelli identity:
\be
\on{rdet}\Big(\sum_{a=1}^M x_{ja}p_{ia}+h(M-i)\delta_{ij}\Big)_{i,j=1}^M=
\det X\ \det P.
\ee
If $h=0$ then all entries of $X$ and $P$ commute and 
the Capelli identity reads $\det(XP)=\det(X)\det(P)$.
Therefore, the Capelli identity can be thought of as a ``quantization'' of the
identity $\det(AB)=\det(A)\det(B)$ for square matrices
$A,B$ with commuting entries.

\subsection{Proof of Theorem \ref{main}}\label{proof}
We denote 
\be
E_{ij,a}:=x_{ja}p_{ia}/(u-z_a).
\ee
We obviously have 
\be
[E_{ij,a},E_{kl,b}]=\delta_{ab}(\delta_{kj}(E_{il,a})'-\delta_{il}(E_{kj,a})'),
\ee
where the prime denotes the formal differentiation with respect to $u$.

Denote also $F^{1}_{jk,a}=-E_{jk,a}$ and $F^{0}_{jj,0}=(p_u-\la_j)$.

Expand $\on{rdet}(G)$. We get an alternating sum of terms,
\bean\label{expand}
\on{rdet}(G_h)=\sum_{\sigma,a,c} (-1)^{\sgn(\sigma)}
F^{c(1)}_{1\sigma(1),a(1)}F^{c(2)}_{2\sigma(2),a(2)}\dots
F^{c(M)}_{M\sigma(M),a(M)},
\eean
where the summation is over all triples $\sigma,a,c$ such that 
$\sigma$ is a permutation of $\{1,\dots,M\}$ and $a,c$ are maps
$a:\ \{1,\dots, M\}\to \{0,1,\dots,N\}$, 
$c:\ \{1,\dots,M\}\to \{0,1\}$ 
satisfying: $c(i)=1$ if $\sigma(i)\neq i$; $a(i)=0$ if and only if $c(i)=0$.

 Let $m$ be a product whose factors are of the form $f(u)$,
   $p_u$, $p_{ij}$, $x_{ij}$ where $f(u)$ are some rational functions in
   $u$. Then the product $m$ will be
   called {\it normally ordered} if all factors of the form $p_u, p_{ij}$
   are on the right from all factors of the form $f(u), x_{ij}$. For
   example,   $(u-1)^{-2} x_{11} p_u p_{11}$ is normally ordered and
   $p_u (u-1)^{-2} x_{11} p_{11}$ is not.

   Given a product $m $ as above, define a new normally ordered
   product $:m:$ as the product of all factors of $m$ in which
   all factors of the form $p_u, p_{ij}$ are placed on the right from all
   factors of the form $f(u), x_{ij}$. For
   example, $ :p_u (u-1)^{-2}x_{11}p_{11}: = (u-1)^{-2}x_{11}p_up_{11}$.

If all variables $p_u$, $p_{ij}$ are moved to the right in the
expansion of $\on{rdet}(G)$ then we get terms obtained by normal
ordering from the terms in \Ref{expand} plus new terms created by the
non-trivial commutators. We show that in fact all new terms cancel
in pairs. 

\begin{lem}\label{no action} For $i=1,\dots,M$, we have
\bean\label{move}
\on{rdet}(G_h)=\sum_{\sigma,a,c} (-1)^{\sgn(\sigma)}
F^{c(1)}_{1\sigma(1),a(1)}\dots 
F^{c(i-1)}_{(i-1)\sigma(i-1),a(i-1)}\Big(:F^{c(i)}_{i\sigma(i),a(i)}\dots 
F^{c(M)}_{M\sigma(M),a(M)}:\Big),
\eean
where the sum is over the same triples $\sigma,a,c$ as in \Ref{expand}.
\end{lem}
\begin{proof}
We prove the lemma by induction on $i$.
For $i=M$ the lemma is a tautology. 
Assume it is proved for $i=M,M-1,\dots,j$, let us prove it for $i=j-1$.

We have 
\bean\label{first}
\lefteqn{F^{1}_{(j-1)r,a} :F^{c(j)}_{j\sigma(j),a(j)}\dots 
F^{c(M)}_{M,\sigma(M),a(M)}:=}\\ 
&&
:F^{1}_{(j-1)r,a}
F^{c(j)}_{j\sigma(j),a(j)}\dots 
F^{c(M)}_{M\sigma(M),a(M)}:+
\sum_{k}:F^{c(j)}_{j\sigma(j),a(j)}\dots 
(-E_{kr,a})'\dots  F^{c(M)}_{M\sigma(M),a(M)}:\ ,\notag
\eean
where the sum is over $k\in\{j,\dots,M\}$ such that $a(k)=a$, 
$\sigma(k)=j-1$ and $c(k)=1$.

We also have 
\bean\label{second}
\lefteqn{F^{0}_{(j-1)(j-1),0}:F^{c(j)}_{j\sigma(j),a(j)}\dots 
F^{c(M)}_{M\sigma(M),a(M)}:=}\\ &&
:F^{0}_{(j-1)(j-1),0}F^{c(j)}_{j\sigma(j),a(j)}\dots 
F^{c(M)}_{M\sigma(M),a(M)}:+\sum_{k}
:F^{c(j)}_{j\sigma(j),a(j)}\dots (-E_{k\sigma(k),a(k)})'\dots  
F^{c(M)}_{M\sigma(M),a(M)}:\ ,\notag
\eean
where the sum is over $k\in\{j,\dots,M\}$ such that $c(k)=1$.

Using \Ref{first}, \Ref{second}, rewrite each term in \Ref{move} with $i=j$.
Then the $k$-th term obtained by using \Ref{first} applied 
to the term labeled by
$\sigma,c,a$ with $c(j-1)=0$ 
cancels with the $k$-th obtained by using \Ref{second} applied to the term labeled by $\tilde
{\sigma},\tilde c, \tilde a$ defined by the following rules. 
\be
\tilde\sigma(i)=\sigma(i) \ \ (i\neq j-1,k), \qquad \tilde\sigma(j-1)=j-1,
 \qquad \tilde \sigma(k)=\sigma(j-1),
\ee
\be
\tilde c(i)=c(i) \ \ (i\neq j-1), \qquad \tilde c(j-1)=0,
\ee
\be
\tilde a(i)=a(i) \ \ (i\neq j-1), \qquad \tilde a(j-1)=0.
\ee
After this cancellation we obtain the statement of the lemma for $i=j-1$.
\end{proof}

\begin{rem} The proof of Lemma \ref{no action} implies that if the matrix 
 $\sigma G_h$ 
 is obtained from $G_h$ by permuting the rows of $G_h$ 
 by a permutation $\sigma$  
  then $\on{rdet}(\sigma G_h)=(-1)^{\sgn(\sigma)}\on{rdet}(G_h)$.
\end{rem}

Consider the linear isomorphism $\phi_h:\ A_h\to A_0$ which sends any 
normally ordered monomial $m$ in $A_h$ to the same monomial $m$ in $A_0$.

By \Ref{move} with $i=1$, the image $\phi_h(\on{rdet}(G_h))$ 
does not depend on $h$ and therefore can be computed at $h=0$. Therefore
Theorem \ref{main} for all $h$ follows from Theorem \ref{main} for $h=0$.
Theorem \ref{main} for $h=0$ follows from formula \Ref{gauss}.

\section{The $(\gl_M,\gl_N)$ duality and the Bethe subalgebras}
\label{bethe sec}
\subsection{Bethe subalgebra}
Let $E_{ij}$, $i,j = 1, \dots, M$, be the standard generators
of $\gl_{M}$. Let $\h$ be the Cartan subalgebra of $\gl_M$,
\be
\h = \oplus_{i=1}^M \C\cdot E_{ii}.
\ee
We denote $U\gl_M$ the universal enveloping algebra of $\gl_M$.

For $\mu \in \h^*$, and a $\gl_M$ module $L$ denote by $L[\mu]$ the
vector subspace of $L$ of vectors of weight
$\mu$,
\bea
L[\mu] &= & \{
v \in L\ {} |\ {} hv\, = \,\langle \mu, h \rangle \,v
\ {\rm for\ any }\ h\in \h \}.
\eea
We always assume that $L=\oplus_\mu L[\mu]$.

For any integral dominant weight
$\La \in \h^*$, denote by $L_\La$ the finite-dimensional 
irreducible $\gl_{M}$-module with
highest weight $\La$.

Recall that we fixed sequences of complex numbers $\bs
z=(z_1,\dots,z_N)$ and $\bs \la=(\la_1,\dots,\la_M)$. From now on we
will assume that $z_i\neq z_j$ and $\la_i\neq \la_j$ if $i\neq j$.

For $i,j=1,\dots,M,$, $a=1,\dots,N$, let
$E_{ji}^{(a)}=1^{\otimes(a-1)}\otimes E_{ji}\otimes
1^{\otimes(N-a)}\in (U\gl_M)^{\otimes N}$.

Define the $M\times M$ matrix $\widetilde G=\widetilde G(M,N,\bs z,\bs \la,u)$ by
\be
\widetilde G(M,N,\bs z,\bs \la,u):=
\Big((\frac{\pr}{\pr u}
-\la_i)\delta_{ij}-\sum_{a=1}^N\frac{E_{ji}^{(a)}}{u-z_a}\Big)_{i,j=1}^M.
\ee 
The entries of $\widetilde G$ are differential operators in $u$ 
whose coefficients 
are rational functions in $u$ with values in $(U\gl_M)^{\otimes N}$. 

Write
\be
\on{rdet}(\widetilde  G(M,N,\bs z,\bs \la,u))=\frac{\pr^{M}}{\pr u^{M}}+
\widetilde G_1(M,N,\bs z,\bs \la,u)\frac{\pr^{M-1}}{\pr u^{M-1}}+\dots+
\widetilde G_{M}(M,N,\bs z,\bs \la,u).
\ee
The coefficients $\widetilde 
G_i(M,N,\bs z,\bs \la,u)$, $i=1,\dots,M$, are called
{\it the transfer matrices of the Gaudin model}.  The transfer
matrices are rational functions in $u$ with values in
$(U\gl_M)^{\otimes N}$.

The transfer matrices commute: 
\be
[\widetilde G_i(M,N,\bs z,\bs \la,u),\widetilde G_j(M,N,\bs z,\bs \la,v)]=0,
\ee
for all $i,j,u,v$, see
\cite{T} and Proposition 7.2 in \cite{MTV}.

The transfer matrices clearly 
commute with the diagonal action  of ${\frak h}$ 
on $(U\gl_M)^{\otimes N}$.

For $i=1,\dots,M$, it is clear that $\widetilde G_i(M,N,\bs z,\bs \la,u)
\prod_{a=1}^N(u-z_a)^i$ is a polynomial in $u$ whose coefficients are
pairwise commuting elements of $(U\gl_M)^{\otimes N}$.  Let 
$\mc G(M,N,\bs z,\bs \la)\subset (U\gl_M)^{\otimes N}$ 
be the commutative subalgebra
generated by the coefficients of polynomials $\widetilde G_i(M,N,\bs z,\bs
\la, u) \prod_{a=1}^N(u-z_a)^i$, $i=1,\dots,M$. 
We call the subalgebra
$\mc G(M,N,\bs z,\bs \la)$ the {\it Bethe subalgebra}.

Let $\mc G(M,\bs \la)\subset U\gl_M[t]$ 
be the subalgebra considered in the introduction.
  Let $U\gl_M[t] \to (U\gl_M)^{\otimes N}$ be the algebra homomorphism
  defined by $E_{ij}\otimes t^n \mapsto \sum_{a=1}^NE_{ij}^{(a)}z_a^n$.
  Then the subalgebra
  $\mc G(M,N,\bs z,\bs \la)$ is the image of the subalgebra $\mc G(M,\bs \la)$
  under that homomorphism.

The Bethe subalgebra clearly acts on any $N$-fold tensor products of 
$\gl_M$ representations.

\medskip

Define the {\it Gaudin Hamiltonians}, 
$H_a(M,N,\bs z,\bs \la)\subset (U\gl_M)^{\otimes N}$, 
$a=1,\dots,N$, by the formula 
\be
H_a(M,N,\bs z,\bs \la)=
\sum_{b=1,b\neq a}^N\frac{\Omega^{(ab)}}{z_a-z_b}+\sum_{b=1}^M 
\la_bE_{bb}^{(a)},
\ee
where $\Omega^{(ab)}:=\sum_{i,j=1}^M E_{ij}^{(a)}E_{ji}^{(b)}$.

Define the {\it dynamical Hamiltonians} 
$H^\vee_a(M,N,\bs z,\bs \la)\subset (U\gl_M)^{\otimes N}$,
$a=1,\dots,M$, by the formula
\bea
H^\vee_a(M,N,\bs z,\bs \la)\,  =\
\sum_{b=1,\ b \ne a}^M
\frac {(\sum_{i=1}^N E_{ab}^{(i)})( \sum_{i=1}^N E_{ba}^{(i)})-\sum_{i=1}^N 
E_{aa}^{(i)}}{\la_a -\la_b}\ 
+ \ \sum_{b=1}^N\ z_b\, E_{aa}^{(b)}\ .
\eea

It is known that the Gaudin Hamiltonians  and the dynamical Hamiltonians 
are in the Bethe subalgebra, see e.g. Appendix B in \cite{MTV}:
\be
H_a(M,N,\bs z,\bs \la)\in \mc G(M,N,\bs z,\bs \la), \qquad H^\vee_b(M,N,\bs z,\bs \la)\in\mc G(M,N,\bs z,\bs \la),
\ee
$a=1,\dots,N$, $b=1,\dots,M$. 

\subsection{The $(\gl_M,\gl_N)$ duality}
Let $\Vb\MMM=\C[x_1,\dots, x_M]$ be the space of polynomials of $M$ variables.
We define the $\gl_M$-action on $\Vb\MMM$ by the formula
\be 
E_{ij}\mapsto x_i\frac \partial{\partial x_j}.
\ee

Then we have an isomorphism of $\gl_M$ modules
\be
\Vb\MMM\,=\,\bigoplus_{m=0}^\infty\,L_{m}\MMM\,
\ee
the submodule $L_{m}\MMM$ being spanned by homogeneous 
polynomials of degree $m$.
The submodule $L_{m}\MMM$ is the irreducible $\gl_M$ module 
with highest weight $(m,0,\dots,0)$ and
highest weight vector  $x_1^{m}$.

Let
$\Vb^{( M,N )} = \C [x_{11},\dots,x_{1N},\dots,x_{M1},\dots,x_{MN}]$
be the space of polynomials of $MN$ commuting variables. 

Let
$
\pi\MMM:\ (U\gl_M)^{\otimes N} \to \on{End}(\Vb^{( M,N )})
$
be the algebra homomorphism defined by 
\be
 E_{ij}^{(a)} \mapsto x_{ia} \frac \partial{\partial x_{ja}}.
\ee
In particular, we define the $\gl_M$ action on $\Vb^{(M,N)}$ by the formula
\be
 E_{ij}
\mapsto \sum_{a=1}^N\,x_{ia} \frac \partial{\partial x_{ja}}.
\ee

Let 
$
\pi\NNN:\ (U\gl_N)^{\otimes M} \to \on{End}(\Vb^{ (M,N)})
$
be the algebra homomorphism defined by
\be
 E_{ij}^{(a)} \mapsto  x_{ai} \frac \partial{\partial x_{aj}}.
\ee
In particular, we define the  $\gl_N$ action on $\Vb^{( M,N )}$ by the formula
\be
E_{ij}
\mapsto \sum_{a=1}^M\,x_{ai} \frac \partial{\partial x_{aj}}.
\ee

We have isomorphisms of algebras, 
\bean\label{Miso}
\bigl(\C[x_1,\dots, x_M]\big)^{\otimes N}\!\!\!\to \Vb^{( M,N )},
&& 
1^{\otimes(j-1)}\otimes x_i \otimes 1^{\otimes(N-j)} \mapsto  x_{ij},\notag
\\
\bigl(\C[x_1,\dots, x_N]\big)^{\otimes M}\!\!\!\to \Vb^{( M,N )},
&& 
1^{\otimes(i-1)}\otimes x_j \otimes 1^{\otimes(M-i)} \mapsto  x_{ij}.
\eean
Under these isomorphisms the space 
$\Vb^{(M,N)}$ is isomorphic to $(\Vb\MMM) ^{\otimes N}$ as a $\gl_M$ module
and  to $({\Vb}^{( N )})^{\otimes M}$ as a $\gl_N$ module.

Fix $\bs n=(n_1,\dots,n_N)\in\Z^{N}_{\geq 0}$ and $\bs m=(m_1,\dots,
m_M)\in\Z_{\geq 0}^{M}$ with $\sum_{i=1}^N n_i = \sum_{a=1}^M m_a$.
The sequences $\bs n$ and $\bs m$ naturally 
correspond to integral $\gl_N$ and $\gl_M$
weights, respectively.

Let $\bs L_{\bs m}$ and  $\bs L_{\bs n}$ be $\gl_N$ and $\gl_M$ modules, 
respectively, defined by the formulas
\be 
\bs L_{\bs m}=\otimes_{a=1}^M L_{m_a}\NNN, \qquad \bs L_{\bs n}=\otimes_{b=1}^N L_{n_b}\MMM.
\ee
The isomorphisms \Ref{Miso} induce an isomorphism of the weight subspaces,
\bean\label{duality isom}
 {\bs L}_{\bs n}[\bs m]\ \simeq  {\bs L}_{\bs m}[\bs n].
\eean

Under the isomorphism \Ref{duality isom} the Gaudin and dynamical 
Hamiltonians interchange,
\bea
\pi\MMM H_{a}(M,N,\bs z,\bs \la)= \pi\NNN H^\vee_{a}(N,M,\bs \la,\bs z)\,, 
 \qquad 
\pi\MMM H^\vee_{b}(M,N,\bs z,\bs \la)=\pi\NNN H_{b}(N,M,\bs \la,\bs z)\,, 
\eea
for $a=1,\dots, N$,\  $b=1,\dots,M$, see \cite{TV}.

We prove a stronger statement that 
the images of $\gl_M$ and $\gl_N$ Bethe subalgebras in 
$\on{End}(\Vb^{ (M,N)})$ are the same.
\begin{theorem}\label{dual}
We have 
\be
\pi\MMM (\mc G(M,N,\bs z,\bs\la))=\pi\NNN (\mc G(N,M,\bs\la,\bs z)).
\ee
Moreover, we have
\bea
\prod_{a=1}^N (u-z_a) \pi\MMM \on{rdet}(\widetilde G(M,N,\bs z,\bs \la,u))=
\sum_{a=1}^N\sum_{b=1}^M A_{ab}\MMM\ u^a \frac{\pr^b}{\pr u^b}, \\ 
\prod_{b=1}^M (v-\la_b) \pi\NNN \on{rdet}(\widetilde G(N,M, \bs \la,\bs z,v))=
\sum_{a=1}^N\sum_{b=1}^M A_{ab}\NNN\ v^b \frac{\pr^a}{\pr v^a},
\eea
where $A_{ab}\MMM$, $A_{ab}\NNN$ 
are linear operators independent on $u,v, \pr/\pr u, \pr/\pr v$ and
\be A_{ab}\MMM=A_{ab}\NNN .\ee

\end{theorem}
\begin{proof}
We obviously have 
\bea
\pi\MMM ( \widetilde G(M,N,\bs z,\bs \la,u))=\bar i_{h=1}(G_{h=1}), \\
\pi\NNN  (\widetilde G(N,M,\bs \la,\bs z,v))=\bar i_{h=1}(H_{h=1}),
\eea
where $G_{h=1}$ and $H_{h=1}$ are matrices defined in \Ref{G} and \Ref{H}.

Then the coefficients of the differential operators $\prod_{a=1}^N
(u-z_a)\pi\MMM \on{rdet}( \widetilde G(M,N,\bs z,\bs \la,u))$ and
$\prod_{b=1}^M (v-\la_b) \pi\NNN \on{rdet}(\widetilde G(N,M,\bs
\la,\bs z,v))$ are polynomials in $u$ and $v$ of degrees $N$ and $M$,
respectively, by Theorem \ref{main}. The rest of the theorem follows
directly from Corollary \ref{duality rel}.
\end{proof}

\subsection{Scalar differential operators}
Let $w\in{\bs L}_{\bs n}[\bs m] $ be a common eigenvector of the
Bethe subalgebra $\mc G(M,N,\bs z,\bs\la)$. Then the operator 
$\on{rdet}(\widetilde G(M,N,\bs z,\bs \la,u))$ acting on $w$ 
defines a monic scalar differential operator of order $M$ with rational 
coefficients in variable $u$. Namely, let $D_w(M, N,\bs \la,\bs z)$ 
be the differential operator
given by
\be
D_w(M, N,\bs z, \bs \la,u) = \frac{\pr^{M}}{\pr u^{M}}+
\widetilde G_1^w(M,N,\bs z,\bs \la,u)\frac{\pr^{M-1}}{\pr u^{M-1}}+\dots+
\widetilde G_{M}^w(M,N,\bs z,\bs \la,u), 
\ee
where $\widetilde G_i^w(M,N,\bs z,\bs \la,u)$ is the eigenvalue of the $i$th 
transfer matrix acting on the vector $w$:
\be
\widetilde G_i(M,N,\bs z,\bs \la,u)w=\widetilde 
G_i^w(M,N,\bs z,\bs \la,u)w.
\ee

Using isomorphism  \Ref{duality isom}, consider $w$ as a vector in 
${\bs L}_{\bs m}[\bs n] $. Then 
by Theorem \ref{dual}, $w$ is also a common eigenvector for algebra 
$\mc G(N,M,\bs \la,\bs z)$.
Thus, similarly, the operator 
$\on{rdet}(\widetilde G(N,M,\bs \la,\bs z,v))$ acting on $w$ 
defines a monic scalar differential operator of order $N$,
$D_w(N,M,\bs \la, \bs z,v)$. 

\begin{cor}\label{eigen dual}
We have
\bea
\prod_{a=1}^N (u-z_a) D_w(M,N,\bs z,\bs \la,u)=
\sum_{a=1}^N\sum_{b=1}^M A^{(M)}_{ab,w}\ u^a \frac{\pr^b}{\pr u^b}, \\ 
\prod_{b=1}^M (v-\la_b) D_w(N,M,\bs \la,\bs z,v)=
\sum_{a=1}^N\sum_{b=1}^M A^{(N)}_{ab,w}\  v^b \frac{\pr^a}{\pr v^a},
\eea
 where $A_{ab,w}^{(M)}$, $A_{ab,w}^{(N)}$ are numbers independent on 
$u,v, \pr/\pr u, \pr/\pr v$. Moreover,
\be
A_{ab,w}^{(M)}=A_{ab,w}^{(N)}.
\ee
\end{cor}
\begin{proof}
The corollary follows directly from Theorem \ref{dual}.
\end{proof}

\medskip

Corollary \ref{eigen dual} was essentially conjectured in Conjecture 5.1 
in \cite{MTV2}.

\begin{rem}
 The operators $D_w(M,N,\bs z,\bs \la)$ are useful objects,
 see \cite{MV1}, \cite{MTV2}, \cite{MTV3}. They have the following three
 properties.

 \begin{enumerate}

 \item[(i)]
 The kernel of $D_w(M,N,\bs z,\bs \la)$ is spanned by the functions
 $p_i^w(u)e^{\la_iu}$, $i=1,\dots,M$, where $p_i^w(u)$ is a polynomial in
 $u$ of degree $m_i$.

 \item[(ii)]
 All finite singular points of $D_w(M,N,\bs z,\bs \la)$
 are $z_1,\dots,z_N$.

 \item[(iii)]
 Each singular point $z_i$ is regular and the
 exponents of $D_w(M,N,\bs z,\bs \la)$ 
 at $z_i$ are  $0,n_i+1,n_i+2,\dots,n_i+M-1$.
 \end{enumerate}

 A converse statement is also true. Namely,
 if a linear differential operator of order $M$ has
 properties (i-iii), then the operator has the form $D_w(M,N,\bs z,\bs
 \la)$ for a suitable eigenvector $w$ of the Bethe subalgebra. This statement 
 may be deduced from Proposition \ref{simple} below.

 We will discuss the properties of such differential operators in
 \cite{MTV4}, cf. also \cite{MTV2} and Appendix A in \cite{MTV3}.
\end{rem}

\subsection{The simple joint spectrum of the Bethe subalgebra}
It is proved in \cite{R}, that 
for any tensor product of irreducible $\gl_M$ modules and for 
generic $\bs z,\bs \la$ the Bethe subalgebra has a simple joint spectrum. 
We give here a proof of this fact in the special case of the tensor product 
$\bs L_{\bs n}$.

\begin{prop}\label{simple}
For generic values of $\bs \la$, the 
joint spectrum of the Bethe subalgebra 
$\mc G(M,N,\bs z,\bs\la)$ acting in ${\bs L}_{\bs n}[\bs m]$ is simple.
\end{prop}
\begin{proof}
We claim that for generic values of $\bs \la$, the 
joint spectrum of the Gaudin Hamiltonians $H_a(M,N,\bs z,\bs \la)$, 
$a=1,\dots, N$, acting in ${\bs L}_{\bs n}[\bs m]$ 
is simple. Indeed fix $\bs z$
and 
consider $\bs \la$ such that 
$\la_1\gg\la_2\gg\dots \gg \la_M\gg 0$. 
Then the eigenvectors 
of the Gaudin Hamiltonians in ${\bs L}_{\bs n}[\bs m]$ 
will have the form $v_1\otimes\dots\otimes v_N+o(1)$, where 
$v_i\in L_{n_i}[\bs m^{( i)}]$ and $\bs m=\sum_{i=1}^N \bs m^{(i)}$. 
The corresponding eigenvalue of $H_{a}(M,N,\bs z,\bs \la)$ will be 
$\sum_{j=1}^M \la_j m_j^{(a)}+O(1)$.

The weight spaces $L_{n_i}\MMM[\bs m_{\bs i}]$ all have 
dimension at most 1 and 
therefore the joint spectrum is simple in this asymptotic zone of parameters. 
Therefore it is simple for generic values of $\bs \la$.
\end{proof}

\end{document}